\documentclass{article} 

\usepackage{amsmath,amsthm}     
\usepackage{graphicx}     
\usepackage{hyperref} 
\usepackage{url}
\usepackage{amsfonts}

\theoremstyle{theorem}
\newtheorem{theorem}{Theorem}

\newtheorem{lem}{Lemma}

\theoremstyle{definition}

\allowdisplaybreaks

\makeatletter
\@addtoreset{footnote}{page}
\makeatother

\begin{document}

\title{Chebyshev Polynomials, Sliding Columns, and the $k$-step Fibonacci Numbers.}

\author{Greg Dresden\\           
Lexington, VA 24450\\            
dresdeng@wlu.edu}                

\maketitle

Let's begin 
by introducing
the Chebyshev polynomials of the first kind, defined here as
$T_0(x) = 1$, 
$T_1(x) = x$, and
$T_n(x) = 2x T_{n-1}(x) - T_{n-2}(x)$ for $n \geq 2$. 
To quote from Benjamin and Walton \cite{BenWal},
\begin{quote}
{\em It's hard to avoid Chebyshev polynomials. They appear in
just about every branch of mathematics, including geometry,
combinatorics, number theory, differential equations,
approximation theory, numerical analysis, and statistics.}
\end{quote}
We will not go into further detail as to their ubiquity; instead,
this article will
illustrate an interesting connection between these 
Chebyshev polynomials and 
the generalized Fibonacci numbers. 
 
The Fibonacci numbers, of course, are even more  well-known
than   the
Chebyshev polynomials; these numbers are defined by 
$F_n = 0$ for $n \leq 0$, $F_1 = 1$, and thereafter each $F_n$ is the sum of
the two previous numbers. A common
generalization is to define the $k$-step Fibonaccis $F_n^{(k)}$ for $k \geq 1$
as follows:
\[
F_n^{(k)} = \begin{cases}
	0, & \text{ if $n \leq 0$}; \\
	1, & \text{ if $n = 1$}; \\
	F_{n-1}^{(k)} + F_{n-2}^{(k)} + \cdots + F_{n-k}^{(k)}, & \text{ if $n \geq 2$.}
	\end{cases}
\]
These are also
called the $k$-generalized Fibonaccis or the $k$-bonaccis, and 
for $k=2$ they give our familiar Fibonaccis: $1, 1, 2, 3, 5, 8, 13, \dots$. 
For $k=3$ they produce what are commonly called 
the Tribonaccis, and for $k=4$ the 
Tetranaccis, and so on.

Returning now to our Chebyshev polynomials, we write out the first seven of them in  Table \ref{tab1}.
	\begin{table}[h]	
	\begin{center}
	\begin{tabular}{rclllllll}
	\hline
	$T_0(x)$ & = & $1$ & & & & & & \\
	$T_1(x)$ & = & $1x$ & & & & & & \\
	$T_2(x)$ & = & $2x^2$ &$-1$  & & & & & \\
	$T_3(x)$ & = & $4x^3$ &$-3x$  & & & & & \\
	$T_4(x)$ & = & $8x^4$ &$-8x^2$  &$+1$ & & & & \\
	$T_5(x)$ & = & $16x^5$ &$-20x^3$  &$+5x$ & & & & \\
	 \hline
	\end{tabular} 
	\end{center}
	\caption{The Chebyshev polynomials $T_n(x)$}\label{tab1}
	\end{table}
It's interesting to note that the sum of the coefficients along the rising diagonals
of Table \ref{tab1}
gives the  Fibonacci numbers;
Bergum, Wagner, and Hoggatt proved this in 1975 
\cite{BWH}.
As it turns out, if we take
the sum of coefficients along rising diagonals of steeper and steeper slopes, we get
the Tribonaccis, the Tetranaccis, and so on. Here is the precise statement:

\begin{theorem}\label{the1}
For the table of Chebyshev polynomials as given above, and for $k\geq 1$, 
the sums of coefficients along the lines of
slope $k-1$ equal  the $k$-step Fibonaccis.
\end{theorem}

A few comments are in order. First, this ``sum along rising diagonals" procedure might remind the reader of how the Fibonaccis can also be found by taking sums 
along the (gently) rising diagonals of Pascal's triangle. This is no 
coincidence; both the Fibonacci numbers and the Chebyshev polynomials
can be written in terms of binomial coefficients. And second, 
Paolo Serafini in an unpublished note \cite{Ser} from 2013 used generating functions to come up with the
following summation formula for the generalized Fibonacci numbers (with, yes, binomial coefficients),
\begin{equation*}
F_n^{(k)} = \sum_{r \leq \lfloor (n-1)/(k+1)\rfloor} (-1)^r 2^{n-2-(k+1)r} 
     \left( \binom{n-1-kr}{r} + \binom{n-2-kr}{r-1} \right),
\end{equation*}
and though Serafini apparently did not realize it at the time, those summands 
are exactly the coefficients of the $x^{n-1-(k+1)r}$ term in 
the Chebyshev polynomial $T_{n-1-(k-1)r}$ and so with a bit of 
effort one can tease out that Serafini's expression  is indeed the sum of rising diagonals of slope $k-1$ from
Table \ref{tab1}.

Serafini's formula inspired us to find an easier and more 
natural approach. As a result, our proof uses nothing more than a simple sliding argument 
and the following key insight: 
the $k$-step 
Fibonaccis also satisfy $F_n^{(k)} = 2F_{n-1}^{(k)}  - F_{n-1-k}^{(k)}$ which looks suspiciously like the recurrence relation 
$T_n = 2xT_{n-1} - T_{n-2}$
for the Chebyshevs. 

One issue is that there is only one sequence of Chebyshev polynomials but
infinitely many sequences of $k$-step Fibonacci numbers. Let us 
now give each $k$-step Fibonacci sequence $F_n^{(k)}$ its own  $k$-generalized or $k$-step
Chebyshev polynomials, which we call  $T_n^{(k)}$ for $k \geq 1$ and which we define as follows:
\[
\mbox{for $k \geq 1,$\ \ } T_n^{(k)}(x) = \begin{cases}
	0, & \text{ if $n < 0$}; \\
	1, & \text{ if $n = 0$}; \\
	x, & \text{ if $n = 1$}; \\
	2x T_{n-1}^{(k)}(x) - T_{n-1-k}^{(k)}(x), & \text{ if $n \geq 2$.}
	\end{cases}
\]
Note that when $k=1$ we get our standard Chebyshev polynomials as defined
earlier. For reasons to become clear in a moment, we also need to define $0$-step Chebyshevs $T_n^{(0)}$ but with slightly different initial conditions. We set 
$T_0^{(0)}(x)=1$,
$T_1^{(0)}(x)=x-1$, and 
$T_n^{(0)}(x)=   2x T_{n-1}^{(0)}(x) - T_{n-1-0}^{(0)}(x)$ for 
$n>1$.
That last recurrence relation 
follows the general pattern of recurrence relations for $T_n^{(k)}$,
but in our case of $k=0$ it 
simplifies nicely to give us the elegant formula $T_n^{(0)}(x) = (2x-1)^{n-1}(x-1)$. 
If we write out these polynomials $T_n^{(0)}$ in a table we see a striking similarity to 
Table \ref{tab1} earlier;
	\begin{table}[h]	
	\begin{center}\begin{tabular}{rclllllll}
	\hline
	$T_0^{(0)}(x)$ & = & $1$ & & & & & & \\
	$T_1^{(0)}(x)$ & = & $1x$ & $-1$& & & & & \\
	$T_2^{(0)}(x)$ & = & $2x^2$ &$-3x$  &$+1$& & & & \\
	$T_3^{(0)}(x)$ & = & $4x^3$ &$-8x^2$  & $+5x$& $-1$& & & \\
	$T_4^{(0)}(x)$ & = & $8x^4$ &$-20x^3$  &$+18x^2$ & $-7x$&$+1$ & & \\
	$T_5^{(0)}(x)$ & = & $16x^5$ &$-48x^4$  &$ +56x^3$ & $-32x^2$ & $+9x$& $-1$& \\
	 \hline
	\end{tabular}\end{center}
	\caption{The $k$-step Chebyshev polynomials $T_n^{(k)}(x)$ for $k=0$}\label{tab2}
	\end{table}
the columns of terms in Table \ref{tab1} for
$T_{n}(x)$ are exactly the same as the columns 
in Table \ref{tab2} for $T_n^{(0)}(x)$, except that  
each column in Table \ref{tab1} is pushed down compared to the columns in Table \ref{tab2}.
This is always the case, as we describe here:

\begin{lem}\label{lem3} 
Given any $k\geq 0$, the table of $k$-step 
Chebyshev polynomials 
$T_{n}^{(k)}(x)$ can be obtained from Table \ref{tab2} for 
$T_{n}^{(0)}(x)$ by sliding down each column of terms 
in Table \ref{tab2} so that each column now starts $k+1$ terms below the 
column to its left.
\end{lem}
In other words,  Lemma \ref{lem3} says 
that  we can form a table for the $k$-step Chebyshevs $T_{n}^{(k)}(x)$
by starting with  Table \ref{tab2}, leaving the first column alone,  sliding
the second column down  $k$ steps, 
the third column down $2k$ steps, the fourth column down $3k$ steps, and so on. (Hence, the name {\em $k$-step} Chebyshevs.)

\begin{proof}[Proof of Lemma \ref{lem3}]
In what follows, we let $k$ be any non-negative integer. 
Note that by our recurrence relation 
$T_{n}^{(k)}(x) = 2x T_{n-1}^{(k)}(x) - T_{n-1-k}^{(k)}(x)$ and our
initial conditions (for both $k=0$ and $k>0$), then the leading term
of $T_{n}^{(k)}(x)$ is always $2^{n-1}x^n$. And again because of the recurrence
relation, then each term that's not in the first column is a sum of the term one step
above it (times $2x$) and the term $k+1$ steps above it and one step to the left (times $-1$). 
But recall that the columns 
in $T_{n}^{(k)}(x)$ are all pushed $k+1$ steps down from the columns to their left, and so if we slide
those columns back up $k$ terms relative to the previous column, then we are now in the situation where each term
in these columns
is dependent on the term above it (times $2x$) and the term above and to the left (times $-1$).
But this will give us the 
table for $T_{n}^{(0)}(x)$ in which each row is $2x-1$ times the previous row, which means
(since the initial conditions are the same) that the terms in the columns are the same for 
$T_{n}^{(k)}(x)$ as for $T_{n}^{(0)}(x)$, as desired.
\end{proof}

We conclude with our final arguments to prove the main theorem. 

\begin{proof}[Proof of Theorem \ref{the1}]
Recall that the $k$-step Chebyshevs have the recurrence relation 
$T_n^{(k)}(x) = 2x T_{n-1}^{(k)}(x) - T_{n-1-k}^{(k)}(x)$,
which when $x$ is replaced by $1$ becomes:
\[
T_n^{(k)}(1) = 2 T_{n-1}^{(k)}(1) - T_{n-1-k}^{(k)}(1).
\]
Now, the $k$-step Fibonaccis are commonly described as having the recurrence
relation $F_n^{(k)} = F_{n-1}^{(k)} + F_{n-2}^{(k)} + \cdots + F_{n-k}^{(k)}$,
but a simple substitution gives the relation:
\[
F_n^{(k)} = 2F_{n-1}^{(k)} - F_{n-1-k}^{(k)}
\]
Thus, so long as the $F_n^{(k)}$'s and the $T_n^{(k)}(1)$'s have the same initial
conditions (and they do, so long as $k\geq 1$), then they are identical sequences of numbers. 

The final step is to note that 
each number $T_n^{(k)}(1)$ is the sum of the coefficients of 
$T_n^{(k)}(x)$ along a horizontal line. But by Lemma \ref{lem3} these coefficients are the same
as those for $T_n(x) = T_n^{(1)}(x)$, except that each column is slid up ($k-1$ steps for the second column, $2k-2$ for the third column, and so on) such 
that what used to be a horizontal line in $T_n^{(k)}(x)$ is now a line of slope $k-1$ in
$T_n(x)$, and this is exactly what we need to  conclude our proof of
Theorem \ref{the1}.
\end{proof}

(An obvious next step would be to investigate the diagonals of the Chebyshev polynomials of the {\em second} kind, but that will have
to wait for a second paper).

{\small 
\noindent {\footnotesize {\bf GREG DRESDEN} (MR Author ID: 623871) received his
 Ph.D.~from the University of Texas in 1997 and now lives and works in the Blue Ridge Mountains with his wife and two children. As a living kidney donor, he advocates for everyone to learn more about living organ donation. 
}

\end{document}